\documentclass[12pt, halfline, a4paper]{ouparticle}
\usepackage{microtype}
\usepackage{amsmath, amsthm, amssymb}
\usepackage{enumitem}
\usepackage{xcolor}
\usepackage{cleveref}

\renewcommand{\epsilon}{\varepsilon}
\renewcommand{\phi}{\varphi}

\def\dotminus{\mathbin{\ooalign{\hss\raise1ex\hbox{.}\hss\cr\mathsurround=0pt$-$}}}
\def\dotplus{\mathbin{\ooalign{\hss\raise1.5ex\hbox{.}\hss\cr\mathsurround=0pt$+$}}}

\DeclareMathOperator{\diam}{diam}
\DeclareMathOperator{\dist}{dist}
\DeclareMathOperator{\stb}{stb}
\DeclareMathOperator{\Th}{Th}

\newtheorem{theorem}{Theorem}[section]
\newtheorem{corollary}[theorem]{Corollary}
\newtheorem{fact}[theorem]{Fact}
\newtheorem{lemma}[theorem]{Lemma}
\newtheorem{proposition}[theorem]{Proposition}

\theoremstyle{definition}
\newtheorem{definition}[theorem]{Definition}
\newtheorem*{question}{Question}
\newtheorem{remarks}[theorem]{Remarks}

\begin{document}

    \title{A note on $\epsilon$-stability}
    \author{
        \name{Nicolas Chavarria}
        \address{University of Waterloo, Waterloo,\\Ontario, Canada}
        \email{nchavarr@uwaterloo.ca}
    }
    \date{\today}

    \abstract{We study $\epsilon$-stability in continuous logic. We first consider stability in a model, where we obtain a definability of types result with a better approximation than that in the literature. We also prove forking symmetry for $\epsilon$-stability and briefly discuss finitely satisfiable types. We then do a short survey of $\epsilon$-stability in a theory. Finally, we consider the map that takes each formula to its ``degree'' of stability in a given theory and show that it is a seminorm. All of this is done in the context of a first-order formalism that allows predicates to take values in arbitrary compact metric spaces.}

    \maketitle

    \section{Introduction}

    In \cite{chavarriaConantPillay}, Conant, Pillay and I proved a continuous version of Malliaris and Shelah's stable regularity lemma (vid. \cite{malliarisShelah}), which in turn is a version of Szemer{\'e}di's regularity lemma. The objects we study are structures of the form $\langle V,W,f\rangle$, where $V$ and $W$ are \emph{finite} sets, and $f:V\times W\to\mathbb R$ is $k$-$\epsilon$-stable for some $0<k<\omega$ and $\epsilon>0$. We do this by first obtaining a structural result for $\epsilon$-stable formulas in $\aleph_0$-saturated structures, and then using pseudofinite methods to transfer this to the finite context and thus obtain an asymptotic regularity result for these finite bipartite ``weighted'' graphs. 

    In the context of finite structures, it makes little sense to ask a formula $\phi(x,y)$ to be $k$-$\epsilon$-stable for all $\epsilon>0$. Indeed, the following fact is easy to verify, as it is nothing more than an application of the pigeonhole principle.
    \begin{fact}
        Suppose $M$ is a finite (continuous) structure and $\phi(x,y)$ is $k$-$\epsilon$-stable for all $\epsilon>0$. Then, for all $a_0,\ldots,a_{k-1},a_k\in M^x$ and $b_0,\ldots,b_{k-1},b_k\in M^y$, there are distinct $i,j\leq k$ such that  $\phi^M(a_i,b_j)=\phi^M(a_j,b_i)$.
    \end{fact}
    \noindent This is clearly a very strong condition imposed on $\phi(x,y)$. Thus, finite and pseudofinite structures provide a setting where it is more rewarding to study the stability of formulas only approximately and not ``fully''.

    With this as motivation, it is the intent of the present article to present a brief treatment of the theory of $\epsilon$-stable formulas in continuous model theory. Organizationally, this is dealt with in the following way.
    \begin{enumerate}
        \item Section 2 deals with the formalism of continuous logic. While at this point doing a brief presentation of the syntax and semantics of continuous logic has become a hackneyed endeavor, we take this opportunity to propose a way to deal with a logic that allows predicates to take values in arbitrary compact metric spaces. The rest of the article is done in this context.
        \item Section 3 covers the study of $\epsilon$-stability in a model. Included here are a proof of the definabilty of types that obtains a better approximation than that available in the currenty literature, as well as a proof of forking symmetry for the $\epsilon$-stable case. The latter goes through finitely satisfiable types, which we define differently to Ben Yaacov in \cite{benYaacov}, and about which we prove a couple of things. All proofs are laid out in detail in this section.
        \item Section 4 is about stability in a theory. This is here mostly for completion, and for contrast with section 3. It includes a brief description of the Cantor-Bendixson analysis of local type spaces, which is again done in the $\epsilon$-stable case as opposed to the fully stable case. While some of these results were more or less implicit in the literature, we make them explicit here.
        \item Section 5 is a brief section that considers the map that takes each $L$-formula to its ``degree'' of stability in a given theory $T$. We prove that this map is a seminorm continuous with respect to diameter.
    \end{enumerate}

    \subsection*{Notation}

    Notation is fairly standard for continuous model theory. Besides minor details, like using angle brackets instead of parentheses for tuples and sequences, there are three things that are of note.
    \begin{enumerate}
        \item We do not speak of sets of ``conditions'' being satisfiable, but rather of sets of formulas. The understanding being that ``satisfied'' means being equal to $0$, so one can think of these formulas as representing the condition obtained by setting them equal to $0$. (Or, whenever a $\dotminus$ appears in the formula, the condition obtained by replacing it with $\leq$.)
        \item For a formula $\phi(x,y)$, we speak of $\phi(x)$-formulas and $\phi^*(y)$-formulas, instead of $\phi$-formulas and $\phi^*$-formulas, and write $S_{\phi(x)}(M)$ and $S_{\phi^*(y)}(M)$, instead of $S_{\phi}(M)$ and $S_{\phi^*}(M)$. This is done to emphasize the free variables involved.
        \item As a stylistic choice, we use $\chi_\delta(X)$ instead of $\|X\|$ to represent the density character of the topological space $X$.
    \end{enumerate}

    Finally, as has become more and more customary in model theory, we write $x,y,\ldots$ for \emph{tuples} of variables, and write $M^x$ to represent the product of the sorts of $M$ corresponding to the sorts of $x$.
    
    \section{Continuous Model Theory}

    In \cite{benyaacovUsvyatsov}, continuous first-order logic is introduced as a formalism in which predicates take values in the interval $[0,1]$ instead of the discrete boolean value set $\{\top,\bot\}$. The restriction to a particular compact interval forces one to use somewhat awkward connectives, like the half-sum or the truncated sum (also called ``dot plus'', $\dotplus$; ``dot minus'', the truncated subtraction, is integral to the expression of inequalities and thus unavoidable). One can easily work around these issues by having predicates take values in their own language-determined compact interval of $\mathbb R$. In fact, other than quantification becoming somewhat unclear, there seems to be no reason why predicates should not take values in, say, (compact subsets of) $\mathbb C$. But then one can also work around the quantifier problem by restricting quantification to those formulas in which it makes sense, i.e. those which take values in $\mathbb R$. This is not a loss, as ultimately the formulas one is most often interested in quantifying are those that measure the distance between the values of distinct predicates. But then this opens the gates to doing metric-valued first-order logic. So long as the metric space in question contains a copy of the reals, then the aforementioned restriction on quantification still applies in a meaningful way. As the extra effort needed is not significant, we take this approach in this article and choose a fixed and ``universal'' metric value space for our predicates. We use the letter $\mathfrak V$ to symbolize it; if the reader prefers to think of $\mathfrak V$ as $\mathbb R$ or $\mathbb C$, he or she is most welcome to.
    
    The space where our formulas will take values, then, is $\mathfrak V=\ell^\infty(\omega)$, endued with its usual norm as a real Banach space. This particular choice of value space is informed by the following ``universality'' property spaces of the form $\ell^\infty(D)$ satisfy.
    
    \begin{fact}\label{fact:metricEmbeddingBanach}
        Let $X$ be a metric space and $D$ a dense subset of $X$. Then $X$ embeds isometrically into $\ell^\infty(D)$. Moreover, its image is closed if $X$ is complete. In particular, every separable metric space embeds isometrically into $\ell^\infty(\omega)$.
    \end{fact}
    
    A simple embedding can be obtained by fixing a point $x_0\in X$. Then $x\mapsto(d(x,d)-d(d,x_0))_{d\in D}$ is an isometric embedding $X\hookrightarrow\ell^\infty(D)$. In fact, following Corollary 1.1 and Proposition 1.3 in \cite{bessagaPelczynski}, one could further assume that the image of $X$ in $\ell^\infty(D)$ is linearly independent\footnote{Potentially under a different embedding.}. Of course, $\ell^\infty(\omega)$ is not unique with this capacity to embed separable metric spaces (read Urysohn sphere). However, we also have the following easy to prove fact that makes working with it desirable, as we can avail ourselves of its coordinates.

    \begin{fact}
        Suppose $X$ is a metric space and $f:X\to\ell^\infty(\omega)$ is such that $x\mapsto f(x)(n)$ is uniformly continuous with the same modulus of continuity for all $n<\omega$. Then $f$ is uniformly continuous.
    \end{fact}

    As to why separable metric spaces in particular, a first-order logic that aspires to have a model theory akin to that of classical $\{\top,\bot\}$-valued logic, i.e. satisfy the Compactness theorem among others, must restrict its formulas to take values in compact spaces. In our case, we desire the formulas to have values in metric spaces, and any compact metric space is known to be separable. Moreover, we have the following fact.

    \begin{fact}[Dugundji's Extension Theorem \cite{dugundji}]
        Let $X$ be a metrizable space, $A\subseteq X$ closed and $V$ a locally convex topological vector space. If $f:A\to V$ is continuous, then there is a continuous $\widetilde f:X\to V$ such that $\widetilde f|_{A}=f$.
    \end{fact}

    As a Banach space, $\mathfrak V$ is automatically locally convex (one definition of which requires the topology to be given by a family of pseudonorms). Therefore, if $X$ and $Y$ are separable metric spaces with $X$ complete, and $u:X^n\to Y$ is continuous, then this extends to a continuous map $\widetilde u:\mathfrak V^n\to\mathfrak V$ under any particular isometric embeddings of $X$ and $Y$. (One may also invoke the hyperconvexity of $\mathfrak V$ to show that such extensions exist even when $X$ is not complete, vid. \cite{aronszajnPanitchpakdi} Theorem 1 and Theorem 4.)

    We will write $\mathfrak d$ for the metric on $\mathfrak V$, i.e. $\mathfrak d(\mathfrak r,\mathfrak s)=\|\mathfrak r-\mathfrak s\|_\infty$. We will also fix some isometric embedding $\mathbb R\hookrightarrow\mathfrak V$ so we can think of $\mathfrak d$ as a 2-Lipschitz function $\mathfrak V^2\to\mathfrak V$ (where the domain has the max metric).

    The general syntactic and semantic constructions for the $\mathfrak V$-valued first-order logic we study in this article follow that in \cite{benYaacovBerensteinHensonUsvyatsov} and \cite{benyaacovUsvyatsov}, but with a small amount of extra care. Namely, to each relation symbol $P$ of a language $L$, besides its arity and modulus of uniform continuity, we attach a compact subset $K_P$ of $\mathfrak V$. This is where the interpretations of $P$ in $L$-structures are meant to take values. We also go through the inductive construction of $L$-formulas to observe how a value space $K_\phi$ is attached to each, and how quantification is dealt with.
    \begin{enumerate}
        \item If $P$ is a predicate symbol of $L$ and $t_0,\ldots t_{n-1}$ are compatible $L$-terms, then $\phi\equiv P(t_0,\ldots,t_{n-1})$ is an $L$-formula and $K_\phi:=K_P$.
        \item If $\phi_0,\ldots,\phi_{n-1}$ are $L$-formulas and $u$ is (a symbol representing) a continuous function $\mathfrak V^n\to\mathfrak V$, then $\phi\equiv u(\phi_0,\ldots,\phi_{n-1})$ is an $L$-formula and $K_\phi:=u(\prod_{i<n}K_{\phi_i})$.
        \item If $\psi$ is an $L$-formula such that $K_\psi\subseteq\mathbb R^{\geq 0}$, and $x$ is a variable, then $\phi\equiv\sup_x\psi$ and $\phi\equiv\inf_x\psi$ are $L$-formulas and $K_\phi:=K_\psi$.
    \end{enumerate}
    For the semantics of this logic, we once again refer to \cite{benYaacovBerensteinHensonUsvyatsov} and \cite{benyaacovUsvyatsov}, the standard references for $[0,1]$-valued continuous model theory. The only thing worth emphasizing is that, in any $L$-structure $M$, $\phi^M(x)$ will be a function from $M^x$ to the language-determined compact set  $K_\phi\subseteq\mathfrak V$.
    
    Of course, while the language enforces a first restriction on the subset of $\mathfrak V$ where a given $L$-formula may take its values, an $L$-structure or an $L$-theory may further narrow this space. We will use the following notation in this connection. For a given $L$-structure $M$, $L$-theory $T$, and $L$-formula $\phi(x)$,
    \begin{align*}
        \diam_M(\phi)&:=\sup\{\mathfrak d(\phi^M(a),\phi^M(b)):a,b\in M^x\},\\
        \diam_T(\phi)&:=\sup\{\mathfrak d(\phi^M(a),\phi^M(b)):a,b\in M^x, M\models T\}.
    \end{align*}
    It is not difficult to see that $\diam_M(\phi)=\diam_{\Th(M)}(\phi)$.
    
    Our intention here is mostly to develop the local theory of stability, so we will fix an $L$-formula $\phi(x,y)$ for the majority of the article.

    \section{Stability in a Model}

    Throughout this section, we fix an $L$-structure $M$. The following is Definition B.1 in \cite{benyaacovUsvyatsov}.

    \begin{definition}
        Let $\epsilon>0$. We say that $\phi(x,y)$ is \textbf{$\epsilon$-stable in $M$} if there are no sequences $\langle a_i\rangle_{i<\omega}$ in $M^x$ and $\langle b_i\rangle_{i<\omega}$ in $M^y$ such that $\mathfrak d(\phi^M(a_i,b_j),\phi^M(a_j,b_i))\geq\epsilon$ for all $i<j<\omega$.
    \end{definition}

    \begin{remarks}
        As an immediate consequence of the definition, $\phi(x,y)$ is $\epsilon$-stable in $M$ if and only if $\phi^*(y,x)$ is $\epsilon$-stable in $M$, where $\phi^*$ is just $\phi$ but with the roles of the object and parameter variables swapped. It is also worth noting that $\phi(x,y)$ is always $\epsilon$-stable for any $\epsilon>\diam_M(\phi)$.
    \end{remarks}

    In classical logic, one of the most important facts about local stability is definability of types. The following sequence of lemmas and propositions cover the corresponding notion for $\epsilon$-stability. In particular, the next proposition and its proof are akin to Lemma B.2 in \cite{benyaacovUsvyatsov}, but with a better degree of approximation for the definitions obtained ($2\epsilon$ plus change instead of $3\epsilon$) at the expense of a potentially larger Lipschitz constant.

    \begin{proposition}\label{proposition:stableInStructurePreDefinition}
        Suppose $\phi(x,y)$ is $\epsilon$-stable in $M$. Fix $\gamma>\delta>0$ and let $p\in S_{\phi(x)}(M)$. There are $a_0,\ldots,a_{n-1}\in M^x$ such that, if $b,c\in M^y$ satisfy $\mathfrak d(\phi^M(a_i,b),\phi^M(a_i,c))<\delta$ for all $i<n$, then $\mathfrak d(\phi(p,b),\phi(p,c))<2\epsilon+\gamma$.
    \end{proposition}

    \begin{proof}
        Suppose, aiming for a contradiction, that for every $a_0,\ldots,a_{n-1}\in M^x$, there are $b,c\in M^y$ such that $\mathfrak d(\phi^M(a_i,b),\phi^M(a_i,c))<\delta$ for all $i<n$, but $\mathfrak d(\phi(p,b),\phi(p,c))$ $\geq2\epsilon+\gamma$. We will write $\zeta=(\gamma-\delta)/2$.

        Now, suppose we have found $a_0,\ldots,a_{n-1}\in M^x$ and $b_0,\ldots,b_{n-1},c_0,\ldots,c_{n-1}$ in $M^y$ such that
        \begin{enumerate}
            \item $\mathfrak d(\phi^M(a_i,b_j),\phi^M(a_i,c_j))<\delta$ for all $i\leq j<n$,
            \item $\mathfrak d(\phi^M(a_j,b_i),\phi(p,b_i))<\zeta$ and $\mathfrak d(\phi^M(a_j,c_i),\phi(p,c_i))<\zeta$ for all $i<j<n$,
            \item $\mathfrak d(\phi(p,b_i),\phi(p,c_i))\geq2\epsilon+\gamma$ for all $i<n$.
        \end{enumerate}
        To extend these sequences, first take some $a_n\in M^x$ such that
        \begin{align*}
            \mathfrak d(\phi^M(a_n,b_i),\phi(p,b_i))<\zeta\text{ and }\mathfrak d(\phi^M(a_n,c_i),\phi(p,c_i))<\zeta
        \end{align*}
        for all $i<n$. By our assumption, there are $b_n,c_n\in M^y$ such that
        \begin{gather*}
            \mathfrak d(\phi^M(a_i,b_n),\phi^M(a_i,c_n))<\delta,\ i\leq n;\\
            \mathfrak d(\phi(p,b_n),\phi(p,c_n))\geq2\epsilon+\gamma.
        \end{gather*}
        The sequences $a_0,\ldots,a_{n-1},a_n$, $b_0,\ldots,b_{n-1},b_n$ and $c_0,\ldots,c_{n-1},c_n$ satisfy points 1, 2 and 3 again. By induction, we can find $a_i\in M^x$, $b_i,c_i\in M^y$, $i<\omega$, such that they satisfy the same three properties.

        We define the following two subsets of $[\omega]^2$:
        \begin{align*}
            S_b&=\{\{i<j\}:\mathfrak d(\phi^M(a_i,b_j),\phi(p,b_i))\geq\epsilon+\zeta\}\\
            S_c&=\{\{i<j\}:\mathfrak d(\phi^M(a_i,c_j),\phi(p,c_i))\geq\epsilon+\zeta\}
        \end{align*}
        We claim that $[\omega]^2=S_b\cup S_c$. Indeed, suppose there were some $i<j<\omega$ such that both $\mathfrak d(\phi^M(a_i,b_j),\phi(p,b_i))<\epsilon+\zeta$ and $\mathfrak d(\phi^M(a_i,c_j),\phi(p,c_i))<\epsilon+\zeta$. Then, since $\mathfrak d(\phi^M(a_i,b_j),\phi^M(a_i,c_j))<\delta$, we have
        \begin{align*}
            \mathfrak d(\phi(p,b_i),\phi(p,c_i))<2\epsilon+2\zeta+\delta=2\epsilon+\gamma.
        \end{align*}
        This of course contradicts the fact that $\mathfrak d(\phi(p,b_i),\phi(p,c_i))\geq2\epsilon+\gamma$. Therefore, by Ramsey's Theorem, we may assume that either $\mathfrak d(\phi^M(a_i,b_j),\phi(p,b_i))\geq\epsilon+\zeta$ for all $i<j<\omega$, or $\mathfrak d(\phi^M(a_i,c_j),\phi(p,c_i))\geq\epsilon+\zeta$ for all $i<j<\omega$. In the first case, since $\mathfrak d(\phi^M(a_j,b_i),\phi(p,b_i))<\zeta$, $\mathfrak d(\phi^M(a_i,b_j),\phi^M(a_j,b_i))>\epsilon$. In the second case, since $\mathfrak d(\phi^M(a_j,c_i),\phi(p,c_i))<\zeta$, $\mathfrak d(\phi^M(a_i,c_j),\phi^M(a_j,c_i))>\epsilon$. In either case, we get a pair of sequences contradicting $\epsilon$-stability of $\phi(x,y)$.
    \end{proof}

    We can now put the functions $\phi(a_i,y)$---where the $a_i$'s are produced by the proposition above---together to obtain an approximate definition for a given $\phi(x)$-type $p$. This is achieved using the following lemma, which is a quick generalization of Lemma B.3 in \cite{benyaacovUsvyatsov}.

    \begin{lemma}
        Let $X$ be a set, $Y$ a metric space, and $f:X\to\mathbb R$ a bounded function. Suppose $g:X\to Y$ and $\epsilon,\delta>0$ are such that, whenever $x,y\in X$ satisfy $d(g(x),g(y))<\delta$, we have $|f(x)-f(y)|<\epsilon$. Then there is a uniformly continuous function $h:Y\to\mathbb R$ such that
        \begin{align*}
            \sup_{x\in X}|f(x)-h\circ g(x)|\leq\epsilon.
        \end{align*}
        Moreover, $h(Y)\subseteq [\inf f(X),\sup f(X)]$ and it is $D/\delta$-Lipschitz, where $D$ is the diameter of $f(X)$, i.e. $D=\sup f(X)-\inf f(X)$.
    \end{lemma}

    \begin{proof}
        Without loss of generality, we may assume $\inf f(X)=0$ and $\sup f(X)=D$. We may also assume that $Y$ is a normed space (cf. \Cref{fact:metricEmbeddingBanach}). We define $v:Z\to\mathbb R$ by
        \begin{align*}
            v(p)=\text{ choice of point in }\{f(x):x\in X\text{ such that }g(x)=p\},
        \end{align*}
        where $Z=g(X)\subseteq Y$. Note that this is well defined by the choice of domain and the assumptions on $f$. Next, for each $q\in Z$, let
        \begin{align*}
            h_q(p)=v(q)\frac{\dist(p,Y\setminus B(q,\delta))}{\delta},\ p\in Y.
        \end{align*}
        By the triangle inequality, it is clear that $h_q$ is $D/\delta$-Lipschitz. It is also immediately evident that $h_q(p)=0$ for all $p\in Y\setminus B(q,\delta)$. Now, because $Y$ is a normed space, $\dist(p,Y\setminus B(q,\delta))\leq \delta=\dist(q,Y\setminus B(q,\delta))$. Therefore, $h_q(q)=v(q)$ and $0\leq h_q(p)\leq v(q)$ throughout.
        
        We finally define
        \begin{align*}
            h(p)=\sup_{q\in Z}h_{q}(p),\ p\in Y.
        \end{align*}
        It is easy to see that $h$ has the same Lipschitz constant as all the $h_q$. It is also clear that its image is in $[0,D]$. Now, let $x\in X$ and let $p=g(x)\in Z$. If $q\in Z$ is such that $d(q,p)\geq\delta$, then clearly $h_q(p)=0$, so such $q$ do not contribute to $h(p)$. On the other hand, if $d(q,p)<\delta$, then for any $y\in X$ such that $g(y)=q$, we have $d(g(y),g(x))<\delta$ and, therefore, $|f(y)-f(x)|<\epsilon$. It follows that $|v(q)-f(x)|<\epsilon$ and so $h(p)\leq f(x)+\epsilon$. On the other hand, $h(p)\geq v(p)>f(x)-\epsilon$. In conclusion, $|f(x)-h(p)|\leq\epsilon$, which yields
        \begin{align*}
            \sup_{x\in X}|f(x)-h\circ g(x)|\leq\epsilon.
        \end{align*}
    \end{proof}

    \begin{remarks}
        If $X$ is a topological space (resp. a metric space) and $g$ is continuous (resp. uniformly continuous), then this lemma says that $f$ is at distance $\epsilon$ (in the supremum metric) from a continuous (resp. uniformly continuous) function.
    \end{remarks}

    \begin{corollary}\label{corollary:2epsilonTypeDefinability}
        Suppose $\phi(x,y)$ is $\epsilon$-stable in $M$. Fix $\gamma>\delta>0$ and let $p\in S_{\phi(x)}(M)$. Then there is a $\phi^*(y)$-formula $\psi(y)$ such that
        \begin{align*}
            \sup_{b\in M^y}\mathfrak d(\psi^M(b),\phi(p,b))\leq2\epsilon+\gamma.
        \end{align*}
        Moreover, $\psi(y)$ is $\diam_M(\phi)/\delta$-Lipschitz (as a function on the metric space $S_{\phi^*(y)}(M)$).
    \end{corollary}

    \begin{proof}
        Let $a_0,\ldots,a_{n-1}\in M^x$ be as given in \Cref{proposition:stableInStructurePreDefinition}. Now, consider $g:M^y\to\mathfrak V^n$ given by $g(b)=\langle\phi^M(a_i,b)\rangle_{i<n}$, where $\mathfrak V^n$ has the max metric. For each $m<\omega$, let $f_m(b)=\phi(p,b)(m)\in\mathbb R$ for all $b\in M^y$. Then, by the previous lemma, there is a $\diam_M(\phi)/\delta$-Lipschitz map $h_m:\mathfrak V^n\to\mathbb R$ such that
        \begin{align*}
            \sup_{b\in M^y}|f_m(b)-h_m\circ g(b)|\leq 2\epsilon+\gamma.
        \end{align*}
        Moreover, all the $h_m$ map to the same bounded subset of $\mathbb R$ (determined by $K_\phi$). Then $u:\mathfrak V^n\to\mathfrak V$ obtained by putting together the $h_m$'s is also $D/\delta$-Lipschitz. It is clear that
        \begin{align*}
            \psi(y)\equiv u(\phi(a_0,y),\ldots,\phi(a_{n-1},y))
        \end{align*}
        is the desired formula.
    \end{proof}

    The last thing we show for $\epsilon$-stability in a model is the corresponding version of Proposition 7.16 in \cite{benyaacovUsvyatsov}. In fact, the latter follows from \Cref{proposition:eStabilitySymmetry} below when $\phi(x,y)$ is assumed to be $\epsilon$-stable in $M$ for \emph{all} $\epsilon>0$. For this we need the following definition. Against the warning in \cite{benYaacov}, Remark 3.2 (i), we speak of ``finitely satisfiable'' types. The definition we give is different to that of idem, Definition 3.1. However, Fact 3.3 (i) vouches for its correctness in this particular context.

    \begin{definition}
        Let $N\succeq M$ and $p\in S_{\phi(x)}(N)$. We say that $p$ is \textbf{finitely satisfiable in $M$} if, for every $n<\omega$, $b_0,\ldots,b_{n-1}\in N^y$ and $\delta>0$, there is some $a\in M^x$ such that $\mathfrak d(\phi^N(a,b_i),\phi(p,b_i))<\delta$ for all $i<n$.
    \end{definition}

    \begin{proposition}\label{proposition:finitelySatisfiableExtension}
        Let $N\succeq N_0\succeq M$ and let $p\in S_{\phi(x)}(N_0)$ be finitely satisfiable in $M$. Then there exists $q\in S_{\phi(x)}(N)$ that extends $p$ and is also finitely satisfiable in $M$.
    \end{proposition}

    \begin{proof}
        For every $n<\omega$, $b_0,\ldots,b_{n-1}\in N^y$, $\mathfrak r_0,\ldots,\mathfrak r_{n-1}\in K_\phi$ and $\delta>0$, let $\psi_{\overline b,\overline{\mathfrak r},\delta}(x)$ be the following $\phi(x)$-formula over $N$:
        \begin{enumerate}
            \item $\delta\dotminus\max\{\mathfrak d(\phi(x,b_i),\mathfrak r_i):i<n\}$, if, for all $a \in M^x$, $\mathfrak d(\phi^N(a,b_i),\mathfrak r_i)\geq\delta$ for some $i<n$.
            \item 0, otherwise.
        \end{enumerate}
        Now, $\psi_{\overline b,\overline{\mathfrak r},\delta}(x)$ is true of any $a\in M^x$ and, since $p$ is finitely satisfiable in $M$,
        \[
            \{\mathfrak d(\phi(x,b),\phi(p,b)):b\in{N_0}^y\}\cup\{\psi_{\overline b,\overline{\mathfrak r},\delta}(x):\overline b\in(N^y)^n,\mathfrak r\in{K_\phi}^n,\delta>0\}
        \]
        is finitely approximately satisfiable. Let $q\in S_{\phi(x)}(N)$ extend this partial type. Clearly, $q$ extends $p$. Suppose, for a contradiction, that $q$ is not finitely satsifiable in $M$. Then there are $b_0,\ldots,b_{n-1}\in N^y$ and $\delta>0$ such that, for all $a\in M^x$, $\mathfrak d(\phi^N(a,b_i),\phi(q,b_i))\geq\delta$ for some $i<n$. Setting $\mathfrak r_i=\phi(q,b_i)$, this means $\psi_{\overline b,\overline{\mathfrak r},\delta}(x)\equiv\delta\dotminus\max\{\mathfrak d(\phi(x,b_i),\mathfrak r_i):i<n\}$. But then this implies $\psi_{\overline b,\overline{\mathfrak r},\delta}(q)=\delta$, contradicting the fact that, by construction, $\psi_{\overline b,\overline{\mathfrak r},\delta}(q)=0$. Thus, $q$ is finitely satisfiable in $M$.
    \end{proof}

    \begin{lemma}
        Let $N\succeq M$ and $q\in S_{\phi(x)}(N)$ finitely satisfiable in $M$. Suppose $\phi(x,y)$ is $\epsilon$-stable in $M$ and $\psi(y)$ is such that $\mathfrak d(\phi(q,b),\psi^M(b))\leq\delta$ for all $b\in M^y$. Then $\mathfrak d(\phi(q,b),\psi^N(b))\leq\delta+\epsilon$ for all $b\in N^y$.
    \end{lemma}

    \begin{proof}
        Suppose, for a contradiction, that $\mathfrak d(\phi(q,b),\psi^N(b))>\epsilon+\delta$ for some $b\in N^y$. We will construct a pair of sequences $\langle a_i\rangle_{i<\omega}$ in $M^x$ and $\langle b_i\rangle_{i<\omega}$ in $M^y$ contradicting the $\epsilon$-stability of $\phi(x,y)$ in $M$. We do this by induction as follows. For $n<\omega$, suppose we are given $a_0,\ldots,a_{n-1}\in M^x$ and $b_0,\ldots,b_{n-1}\in M^y$ such that
        \begin{enumerate}
            \item $\mathfrak d(\phi(q,b_i),\phi^N(a_i,b))>\epsilon$ for all $i<n$, and
            \item $\mathfrak d(\phi^N(a_i,b_j),\phi^N(a_j,b_i))>\epsilon$ for all $i<j<n$.
        \end{enumerate}
        Since $q$ is finitely satisfiable in $M$, there is $a_n\in M^x$ such that
        \begin{align*}
            \mathfrak d(\phi^N(a_n,b),\psi^N(b))>\epsilon+\delta\text{ and }\mathfrak d(\phi^N(a_n,b_i),\phi^N(a_i,b))>\epsilon    
        \end{align*}
        for all $i<n$. Further, there is $b_n\in M^y$ such that $\mathfrak d(\psi^N(b_n),\psi^N(b))<\gamma$, where $\gamma=\mathfrak d(\phi^N(a_n,b),\psi^N(b))-\epsilon-\delta$, and $\mathfrak d(\phi^N(a_i,b_n),\phi^N(a_n,b_i))>\epsilon$ for all $i<n$. Then, since by hypothesis $\mathfrak d(\phi(q,b_n),\psi^N(b_n))\leq\delta$, we obtain
        \begin{align*}
            \mathfrak d(\phi(q,b_n),\phi^N(a_n,b))&\geq\mathfrak d(\phi^N(a_n,b),\psi^N(b))-\mathfrak d(\psi^N(b_n),\psi^N(b))-\mathfrak d(\phi(q,b_n),\psi^N(b_n))\\
            &>\epsilon+\delta-\delta=\epsilon.
        \end{align*}
        Thus, points 1 and 2 are now true for all $i<n+1$ and the induction can go through. This finishes the proof.
    \end{proof}

    \begin{proposition}\label{proposition:eStabilitySymmetry}
        Suppose $\phi(x,y)$ is $\epsilon$-stable in $M$, $p\in S_{\phi(x)}(M)$ and $q\in S_{\phi^*(y)}(M)$. Let $\psi_p(y)$, $\psi_q(x)$, and $\delta_p,\delta_q>0$ be such that
        \begin{enumerate}
            \item $\psi_p(y)$ is a $\phi^*(y)$-formula and $\mathfrak d(\phi(p,b),\psi_p(b))\leq\delta_p$ for all $b\in M^y$, and
            \item $\psi_q(x)$ is a $\phi(x)$-formula and $\mathfrak d(\phi^*(q,a),\psi_q(a))\leq\delta_q$ for all $a\in M^x$.
        \end{enumerate}
        Then $\mathfrak d(\psi_p(q),\psi_q(p))\leq\delta_p+\delta_q+\epsilon$.
    \end{proposition}

    \begin{proof}
        Let $b\in N^y$ realize $q$ in some $N\succeq M$. Let $p'\in S_{\phi(x)}(N)$ be an extension of $p$ finitely satisfiable in $M$. Take any $\gamma>0$. Then, since $\mathfrak d(\phi(p',b),\psi_p^N(b))\leq\delta_p+\epsilon$, there is some $a\in M^x$ such that $\mathfrak d(\phi^N(a,b),\psi_p^N(b))<\delta_p+\epsilon+\gamma/2$ and $\mathfrak d(\psi_q^N(a),\psi_q(p))<\gamma/2$. Now, note that $\mathfrak d(\phi^N(a,b),\psi_q^N(a))\leq\delta_q$ and $\psi_p^N(b)=\psi_p(q)$, so $\mathfrak d(\psi_p(q),\psi_q(p))<\delta_p+\delta_q+\epsilon+\gamma$. Since this is true for all $\gamma>0$, the proof is done.
    \end{proof}

    \begin{remarks}
        \Cref{corollary:2epsilonTypeDefinability} only guarantees $\delta_p,\delta_q\leq 2\epsilon+\gamma/2$ for any choice of $\gamma>0$. Therefore, in the proposition above we would have $\mathfrak d(\psi_p(q),\psi_q(p))\leq 5\epsilon+\gamma$.
    \end{remarks}

    \section{Stability in a theory}

    We now switch gears to study $\epsilon$-stability in a theory, which yields better approximations and uniformity in the definitions of types. Fix a complete $L$-theory $T$. The following are \cite{benyaacovUsvyatsov}, Lemma 7.2 (ii) (or rather the proof thereof) and Definition 7.1, respectively.

    \begin{definition}
        We say that $\phi(x,y)$ is \textbf{$\epsilon$-stable in $T$} if there do not exist a model $M\models T$, $\mathfrak r,\mathfrak s\in K$, and sequences $\langle a_i\rangle_{i<\omega}$ in $M^x$ and $\langle b_i\rangle_{i<\omega}$ in $M^y$ such that, for all $i<j<\omega$, $\phi^M(a_i,b_j)=\mathfrak r$, $\phi^M(a_j,b_i)=\mathfrak s$, and $\mathfrak d(\mathfrak r,\mathfrak s)\geq\epsilon$.
    \end{definition}

    \begin{proposition}\label{proposition:alternateStabilityDefinition}
        $\phi(x,y)$ is $\epsilon$-stable in $T$ if and only if there do not exist a model $M\models T$ and sequences $\langle a_i\rangle_{i<\omega}$ in $M^x$ and $\langle b_i\rangle_{i<\omega}$ in $M^y$ such that, for all distinct $i,j<\omega$, $\mathfrak d(\phi^M(a_i,b_j),\phi^M(a_j,b_i))\geq\epsilon$.
    \end{proposition}

    \begin{proof}
        One direction is obvious. Suppose we had $M$ and sequences as described in the statement of the proposition. Let $\delta>0$ and $U_0,\ldots,U_{n-1}$ be a cover of $K_\phi$ by open sets of diameter $<\delta$. For each pair $\langle i,j\rangle\in n^2$, we let $$S^\delta_{i,j}=\{\{k<l\}\in[\omega]^2:\phi^M(a_k,b_l)\in U_i\text{ and }\phi^M(a_l,b_k)\in U_j\}.$$
        Then $[\omega]^2=\bigcup_{\langle i,j\rangle\in n^2}S^\delta_{i,j}$, so by Ramsey's Theorem we may assume $$\mathfrak d(\phi^M(a_i,b_j),\phi^M(a_{i'},b_{j'}))<\delta\text{ and }\mathfrak d(\phi^M(a_j,b_i),\phi^M(a_{j'},b_{i'}))<\delta,$$
        for all $i<j<\omega$ and $i'<j'<\omega$. The proposition then follows by Compactness.
    \end{proof}

    One of the virtues of stability in a theory as opposed to stability in a structure---besides the better approximations we shall see below---is that it becomes a finitary property and, hence, expressible in the language.

    \begin{definition}
        Let $0<k<\omega$. We say that $\phi(x,y)$ is \textbf{$k$-$\epsilon$-stable in $T$} if there do not exist $M\models T$, $\langle a_i\rangle_{i<k+1}$ in $M^x$ and $\langle b_i\rangle_{i<k+1}$ in $M^y$ such that $\mathfrak d(\phi^M(a_i,b_j),\phi^M(a_j,b_i))\geq\epsilon$ for all distinct $i,j<k+1$.
    \end{definition}

    \begin{proposition}\label{proposition:stableKStable}
        Suppose $\phi(x,y)$ is $\epsilon$-stable in $T$. Then there is $0<k<\omega$ such that $\phi(x,y)$ is $k$-$\epsilon$-stable in $T$.
    \end{proposition}

    \begin{proof}
        This is just a simple application of Compactness.
    \end{proof}

    Naturally, if $\phi(x,y)$ is $\epsilon$-stable in $T$, then it is $\epsilon'$-stable in $T$ for any $\epsilon'\geq\epsilon$. This means that the set $\{\epsilon>0:\phi(x,y)\text{ is }\epsilon\text{-stable in }T\}$ is either of the form $[\epsilon_0,+\infty)$ or of the form $(\epsilon_0,+\infty)$. The latter is the case.

    \begin{fact}[\cite{chavarriaConantPillay} Lemma 3.2]\label{proposition:stabilityOpen}
        There exists some $\epsilon_0=\epsilon_{\phi(x,y),T}>0$ such that $\phi(x,y)$ is $\epsilon$-stable in $T$ if and only if $\epsilon>\epsilon_0$.
    \end{fact}

    The corresponding definabilty of types result for the current context is the following. The proof is only a very slight modification of the proofs of Lemma 7.4 in \cite{benyaacovUsvyatsov}, sprinkled with some of the proof of Lemma 7.2. As such, we leave it to the reader's imagination.

    \begin{fact}
        Suppose $\phi(x,y)$ is $\epsilon$-stable in $T$ and let $p\in S_{\phi(x)}(M)$ for $M\models T$. Then there are $a_n\in M^x$, $n<\omega$, and $K<\omega$ such that, for any $b\in M^x$, $\mathfrak d(\phi^M(a_n,b),\phi(p,b))\leq\epsilon$ for all but \emph{at most} $K$ many indices $n$.
    \end{fact}

    \begin{corollary}
        Suppose $\phi(x,y)$ is $\epsilon$-stable in $T$ and let $p\in S_{\phi(x)}(M)$ for $M\models T$. Then there is a $\phi^*(y)$-formula $\psi(y)$ such that $\mathfrak d(\phi(p,b),\psi^M(b))\leq\epsilon$ for all $b\in M^y$. Moreover, $\psi$ is 1-Lipschitz (as a map from the metric space $S_{\phi^*(y)}(M)$).
    \end{corollary}

    \begin{remarks}\label{remarks:uniformDefinabilityTypes}
        In fact, the $K$ above is independent of the choice of $p\in S_{\phi(x)}(M)$, so we can reformulate the statement of the previous corollary to say that there is a formula $\psi(y,\overline x)$ (with the appropriate shape) such that, for every $p\in S_{\phi(x)}(M)$, there is $\overline a_p\in M^{\overline x}$ such that $\mathfrak d(\phi(p,b),\psi^M(b,\overline a_p))\leq\epsilon$ for all $b\in M^y$.
    \end{remarks}

    The next lemma hinges on the following observation. Let $M\models T$ and
    \begin{align*}
        \mathfrak V^{\oplus M^y}:=\bigoplus_{b\in M^y}^\infty\mathfrak V=\{f:M^y\to\mathfrak V:\sup_{b\in M^y}\|f(b)\|_\infty<\infty\}
    \end{align*}
    with the $\ell^\infty$ norm. Then, identifying each $p\in S_{\phi(x)}(M)$ with the map $b\in M^y\mapsto\phi(p,b)$, we see that $S_{\phi(x)}(M)\subseteq\mathfrak V^{\oplus M^y}$. (Remember that each map $p:M^y\to\mathfrak V$ takes images in $K_\phi$, a compact subset of $\mathfrak V$.) Moreover, $S_{\phi(x)}(M)$ is a metric subspace of $\mathfrak V^{\oplus M^y}$.

    \begin{lemma}\label{lemma:stableEDensity}
        If $\phi(x,y)$ is $\epsilon$-stable in $T$, and $M\models T$, then $S_{\phi(x)}(M)=\bigcup_{\alpha<\lambda}A_\alpha$, where $\lambda\leq\chi_\delta(M^x)$ and $\diam(A_\alpha)\leq 2\epsilon$.
    \end{lemma}

    \begin{proof}
        Let $0<\epsilon'<\epsilon$ be such that $\phi(x,y)$ is $\epsilon'$-stable in $T$ (vid. \Cref{proposition:stabilityOpen}), and let $\gamma=\epsilon-\epsilon'>0$. Let $\psi(y,\overline x)$ be a formula as given in \Cref{remarks:uniformDefinabilityTypes} and note that this is a uniformly continuous map $M^y\times M^{\overline x}\to K_\psi$. Let $D\subseteq M^{\overline x}$ be a dense set, so $|D|=\chi_\delta(M^{\overline x})=\chi_\delta(M^x)$.  Now, for each $\overline d\in D$, let $A_{\overline d}=\overline B(\psi(y,\overline d),\epsilon)\subseteq\mathfrak V^{\oplus M^y}$. By the uniform continuity of $\psi$, given $p\in S_{\phi(x)}(M)$ and $\overline a_p\in M^{\overline x}$ as in \Cref{remarks:uniformDefinabilityTypes}, there is $\delta_p>0$ such that, if $d(\overline a,\overline a_p)<\delta_p$, then $\sup_{b\in M^y}|\psi^M(b,\overline a)-\psi^M(b,\overline a_p)|\leq\gamma$. By the density of $D$, for every $p\in S_{\phi(x)}(M)$, there is $\overline d\in D$ such that $d(\overline d,\overline a_p)<\delta_p$. Then, for all $b\in M^y$, we have
        \begin{align*}
            |\phi(p,b)-\psi^M(b,\overline d)|&\leq|\phi(p,b)-\psi^M(b,\overline a_p)|+|\psi^M(b,\overline a_p)-\psi^M(b,\overline d)|\leq\epsilon'+\gamma=\epsilon.
        \end{align*}
        Thus, $p\in A_{\overline d}$ and the $A_{\overline d}$ cover $S_{\phi(x)}(M)$. Also note that, clearly, we have $\diam(A_{\overline d})\leq 2\epsilon$. Finally, ordering the $\overline d\in D$ with an ordinal and taking intersections with $S_{\phi(x)}(M)$, we obtain the result.
    \end{proof}

    \begin{remarks}
        The proof of the lemma above implicitly uses the fact that $M^x$ is infinite. If $M^x$ is finite, then $\chi_\delta(M^x)=|M^x|=|S_{\phi(x)}(M)|$.
    \end{remarks}

    The previous lemma was added to relate local $\epsilon$-stability with the Cantor-Bendixson analysis of the corresponding type spaces. The definition, lemma and proposition that follow, unsurprisingly, come from \cite{benyaacovUsvyatsov} section 7, where the latter two are phrased in terms of \emph{full} stability (i.e. $\epsilon$-stability for all $\epsilon>0$ simultaneously). Nevertheless, the approximate versions we present are to some extent implicit in the proofs offered in \cite{benyaacovUsvyatsov}.

    \begin{definition}
        Let $X$ be a topometric space. For each ordinal $\alpha$, we define its \emph{$\alpha$-th $\epsilon$-Cantor-Bendixson derivative} $X^{(\alpha)}_\epsilon$ by induction as follows:
        \begin{enumerate}
            \item $X^{(0)}_\epsilon=X$;
            \item $X^{(\alpha+1)}_\epsilon=\displaystyle\bigcap\{F\subseteq X^{(\alpha)}:F\text{ is closed and }\diam(X^{(\alpha)}\setminus F)\leq\epsilon\}$;
            \item $X^{(\alpha)}_\epsilon=\displaystyle\bigcap_{\beta<\alpha}X^{(\beta)}_\epsilon$, for $\alpha$ limit.
        \end{enumerate}
        For $p\in X$, we say that its \emph{$\epsilon$-Cantor-Bendixson rank} is $\alpha$, if $p\in X^{(\alpha)}_\epsilon\setminus X^{(\alpha+1)}_\epsilon$. If no such $\alpha$ exists, then we say that the $\epsilon$-Cantor-Bendixson rank of $p$ is infinite. Finally, we say that $X$ is \emph{$\epsilon$-Cantor-Bendixson analyzable} if the $\epsilon$-Cantor-Bendixson ranks of its points attain a maximum. We often shorten Cantor-Bendixson to CB.
    \end{definition}

    \begin{lemma}
        Suppose $S_{\phi(x)}(M)$ is not $\epsilon$-CB-analyzable for some $\epsilon>0$. Then, for every $\eta\in 2^{<\omega}$, there are $b_\eta\in M^y$ and $\mathfrak r_\eta,\mathfrak s_\eta\in K_\phi$ such that $\mathfrak d(\mathfrak r_\eta,\mathfrak s_\eta)>\epsilon$ and, for all $\sigma\in 2^\omega$,
        \begin{align*}
            \Phi_\sigma(x)=\{\mathfrak d(\phi(x,b_\eta),\mathfrak r_\eta)&\dotminus\gamma_\eta:\eta\subseteq\sigma,\sigma(|\eta|)=0\}\\
            &\cup\{\mathfrak d(\phi(x,b_\eta),\mathfrak s_\eta)\dotminus\gamma_\eta:\eta\subseteq\sigma,\sigma(|\eta|)=1\},
        \end{align*} 
        where $\gamma_\eta=\dfrac{\mathfrak d(\mathfrak r_\eta,\mathfrak s_\eta)-\epsilon}{2}$ and $|\eta|$ is the length of $\eta$, is consistent.
    \end{lemma}

    \begin{proof}
        Let $X$ be the set of points of infinite $\epsilon$-CB-rank in $S_{\phi(x)}(M)$. By hypothesis, this set is non-empty; by construction, it is such that, if $U\subseteq X$ is relatively open, then $\diam(U)>\epsilon$. Thus, we can take $p_0,p_1\in X$ such that $d(p_0,p_1)>\epsilon$. Then, by definition, there is some $b_\emptyset\in M^y$ such that $\mathfrak d(\phi(p_0,b_\emptyset),\phi(p_1,b_\emptyset))>\epsilon$. Let $\mathfrak r_\emptyset=\phi(p_0,b_\emptyset)$ and $\mathfrak s_\emptyset=\phi(p_1,b_\emptyset)$. Now, consider the (``basic'') open subsets
        \begin{align*}
            U_0=\{p\in S_{\phi(x)}(M):\mathfrak d(\phi(p,b_\emptyset),\mathfrak r_\emptyset)<\gamma_\eta\},\\
            U_1=\{p\in S_{\phi(x)}(M):\mathfrak d(\phi(p,b_\emptyset),\mathfrak s_\emptyset)<\gamma_\eta\}.
        \end{align*}
        Their relativizations to $X$ are not empty, as $p_0$ is in the former and $p_1$ in the latter. Thus, they have diameter $>\epsilon$ and so we can find $p_{\langle 0,0\rangle}$, $p_{\langle 0,1\rangle}$ in $U_0\cap X$ and $p_{\langle 1,0\rangle}$, $p_{\langle 1,1\rangle}$ in $U_1\cap X$ such that $d(p_{\langle 0,0\rangle},p_{\langle 0,1\rangle})>\epsilon$ and $d(p_{\langle 1,0\rangle},p_{\langle 1,1\rangle})>\epsilon$. We hence define $b_0,\mathfrak r_0,\mathfrak s_0$ and $b_1,\mathfrak r_1,\mathfrak s_1$ as above. Finally, we continue by induction on the length of sequences. The consistency of $\Phi_\sigma(x)$ then follows from the fact that it is finitely satisfiable, as we can see using the types $p_\eta$ for $\emptyset\neq\eta\subseteq\sigma$.
    \end{proof}

    \begin{proposition}
        If $\phi(x,y)$ is $\epsilon$-stable in $M$, then $S_{\phi(x)}(M)$ is $2\epsilon$-CB-analyzable.
    \end{proposition}

    \begin{proof}
        Suppose not. Let $b_\eta$, $\mathfrak r_\eta$ and $\mathfrak s_\eta$, $\eta\in 2^{<\omega}$, be as in the previous lemma. Now, restrict $L$ to a countable language $L'$ that contains all the symbols for $\phi$. We also take $N\preceq M|_{L'}$ containing all the $b_\eta$ and such that $\chi_\delta(N^x)\leq\aleph_0$. It is clear that $\phi(x,y)$ is $\epsilon$-stable in $N$. For $\sigma\in 2^\omega$, let $p_\sigma\in S_{\phi(x)}(N)$ extend $\Phi_\sigma(x)$. Then clearly, if $\sigma\neq\tau$, we have $d(p_\sigma,p_\tau)>2\epsilon$. Hence, a set of diameter $\leq 2\epsilon$ contains at most one $p_\sigma$. This contradicts \Cref{lemma:stableEDensity}.
    \end{proof}

    \section{Global stability}

    We now turn our attention to all partitioned $L$-formulas simultaneously. We still fix a complete $L$-theory $T$.

    \begin{definition}
        Let $\phi(x,y)$ be an $L$-formula. We define
        \begin{align*}
            \stb_T(\phi(x,y))=\inf\{\epsilon>0:\phi(x,y)\text{ is }\epsilon\text{-stable in }T\}\in\mathbb R^{\geq 0}.
        \end{align*}
    \end{definition}

    \begin{remarks}
        Note that, if $x'\supseteq x$ and $y'\supseteq y$, then $\stb_T(\phi(x',y'))=\stb_T(\phi(x,y))$.
    \end{remarks}

    The next lemma says that $\stb_T$ acts as a seminorm that is continuous with respect to the diameter seminorm. We do not make the ambient space explicit.

    \begin{lemma}
        For all $L$-formulas $\phi(x,y)$ and $\psi(x,y)$, $r\in\mathbb R$ and $\delta>0$,
        \begin{enumerate}
            \item $\stb_T(\phi(x,y))=0$ if $\phi$ is constant in $T$;
            \item $\stb_T(r\phi(x,y))=|r|\stb_T(\phi(x,y))$;
            \item $\stb_T((\phi+\psi)(x,y))\leq\stb_T(\phi(x,y))+\stb_T(\psi(x,y))$;
            \item $\stb_T(\phi(x,y))\leq\diam_T(\phi)$.
        \end{enumerate}
    \end{lemma}

    \begin{proof}
        1. and 4. are obvious and 2. is easy to see.

        Toward 3. let $M\models T$, $\langle a_i\rangle_{i<\omega}$ and $\langle b_i\rangle_{i<\omega}$ sequences in $M^x$ and $M^y$, respectively. Let $\epsilon>0$ and consider the following subsets of $[\omega]^2$.
        \begin{align*}
            A_0=\{\{i<j\}:\mathfrak d(\phi^M(a_i,b_j),\phi^M(a_j,b_i))<\stb_T(\phi(x,y))+\epsilon/2\},\\
            A_1=\{\{i<j\}:\mathfrak d(\phi^M(a_i,b_j),\phi^M(a_j,b_i))\geq\stb_T(\phi(x,y))+\epsilon/2\}.
        \end{align*}
        It is clear that $[\omega]^2=A_0\cup A_1$, so by Ramsey's Theorem there is some $k<2$ and a strictly increasing map $\iota:\omega\to\omega$ such that $\{\iota(i),\iota(j)\}\in A_k$ for any distinct $i,j<\omega$. But note that $k$ cannot be 1, for otherwise $\langle a_{\iota(i)}\rangle_{i<\omega}$ and $\langle b_{\iota(i)}\rangle_{i<\omega}$ would witness $\phi(x,y)$ \emph{not} being $(\stb_T(\phi(x,y))+\epsilon/2)$-stable in $T$, which is impossible. It follows that we may, without loss of generality, assume that $\mathfrak d(\phi^M(a_i,b_j),\phi^M(a_j,b_i))<\stb_T(\phi(x,y))+\epsilon/2$ for \emph{all} distinct $i,j<\omega$. Then, since $\psi(x,y)$ is $(\stb_T(\psi(x,y))+\epsilon/2)$-stable, there are some distinct $i,j<\omega$ such that $\mathfrak d(\psi^M(a_i,b_j),\psi^M(a_j,b_i))<\stb_T(\psi(x,y))+\epsilon/2$. This implies
        \begin{align*}
            \mathfrak d((\phi+\psi)^M(a_i,b_j),(\phi+\psi)^M(a_j,b_i))<\stb_T(\phi(x,y))+\stb_T(\psi(x,y))+\epsilon.
        \end{align*}
        From here, $\stb_T((\phi+\psi)(x,y))\leq\stb_T(\phi(x,y))+\stb_T(\psi(x,y))+\epsilon$. Since $\epsilon$ was arbitrary, the result desired follows.
    \end{proof}

    \begin{definition}
        Let $\mathcal E\in[0,1]$. We say that $T$ is \textbf{$\mathcal E$-stable} if $\stb_T(\phi(x,y))\leq\mathcal E\diam_T(\phi)$ for all partitioned $L$-formulas $\phi(x,y)$. We also write
        \begin{align*}
            \stb(T):=\min\{\mathcal E\in[0,1]:T\text{ is }\mathcal E\text{-stable}\}.
        \end{align*}
        Finally, we say that $T$ is \textbf{stable} if $\stb(T)=0$.
    \end{definition}

    \begin{question}
        Are there any properties of $\stb_T$ as a seminorm that weigh on the structure of the models of $T$?
    \end{question}

    \bibliographystyle{abbrv}
    \bibliography{resources}

\end{document}